\documentstyle[12pt]{article}
\author{Alaoui Youssef}
\title{On the Levi problem with singularities}
\date{ }

\newtheorem{th}{theorem}
\newtheorem{lm}{lemma}
\newcommand{\complexes}{\mbox{I}\!\!\!\mbox{C}}

\begin{document}
\maketitle
\setcounter{page}{1}
{\large 1. Introduction}\\
\\
\hspace*{.1in}Is a complex space X which is the union of an increasing sequence 
\\$X_{1}\subset X_{2}\subset X_{3}\subset \cdots$
of open Stein subspaces itself a Stein space ?\\
\hspace*{.1in}From the begining this question has held great interest in Stein theory.\\
\hspace*{.1in}The special case when $\{X_{j}\}_{j\geq 1}$ is a sequence of Stein
domains in $\complexes^{n}$ had been proved long time ago by Behnke and Stein 
[2].\\
\hspace*{.1in}In 1956, Stein [13] answered positively the question under the 
additional hypothesis that X is reduced and every pair
$(X_{\nu+1},X_{\nu})$ is Runge.\\
\hspace*{.1in}In the general case X is not necessarily holomorphically-convex.\\ 
Fornaess [7],
gave a 3-dimensional example of such situation.\\
\hspace*{.1in}In 1977, Markoe [10] proved the following:\\
Let X be a reduced complex space which the union of an increasing sequence 
$X_{1}\subset X_{2}\subset\cdots\subset X_{n}\subset\cdots$ of Stein domains.\\
\hspace*{.1in}Then X is Stein if and only if $H^{1}(X,O_{X})=0.$\\
\hspace*{.1in}M. Coltoiu has shown in [3] that if $D_{1}\subset 
D_{2}\subset\cdots\subset D_{n}\subset\cdots$ \\is an increasing sequence of Stein
domains in a normal Stein space X, \\then $D=\displaystyle\bigcup_{j\geq 1}D_{j}$ is 
a domain of holomorphy.
(i.e. for each $x\in \partial{D}$ there is $f\in O(D)$ which is not holomorphically 
extendable through x).\\
\hspace*{.1in}The aim of this paper is to prove the following theorems:\\
\begin{th}{-Let X be a Stein normal space of dimension n and
$D\subset \subset X$ an open subset which the union of an increasing
sequence\\ $D_{1}\subset D_{2}\subset\cdots \subset D_{n}\cdots$
of domains of holomorphy in X. Then D is a domain of holomorphy.}
\end{th}
\begin{th}{-A domain of holomorphy D which is relatively compact in a
2-dimensional normal Stein space X
itself is Stein}
\end{th}
\begin{th}{-Let X be a Stein space of dimensinon n and $D\subset X$ an open 
subspace which is the union
of an increasing sequence\\ $D_{1}\subset D_{2}\subset\cdots \subset 
D_{n}\subset\cdots$ of open
Stein subsets of X. Then D itself is Stein, if X has isolated singularities.}
\end{th}
{\bf 2. Preliminaries}\\
\\
\hspace*{.1in}It should be remarked that the statement of theorem 2 is in general\\ false 
if $dim(X)\geq 3$:\\
\hspace*{.1in}Let $X=\{z\in \complexes^{4}: 
z_{1}^{2}+z_{2}^{2}+z_{3}^{2}+z_{4}^{2}=0\},
H=\{z\in \complexes^{4};z_{1}=iz_{2}, z_{3}=iz_{4}\}, \\U=\{z\in X: |z|<1\},$
and $D=U-U\cap H$.\\
X is a Stein normal space of dimension 3 with the singularity only at\\ the
origine. Since D is the complement of a hypersurface on the Stein space U, then
D is
a domain of holomorphy. But D is not Stein. (See [9]).\\
\\
\hspace*{.1in}Let X be a connected n-dimensional, Stein normal space and 
Y be \\the singular locus of X.\\
\hspace*{.1in}There exist finitely many holomorphic maps\\
$$\phi_{j}:X\longrightarrow \complexes^{n}, j=1,\cdots,l$$
with discrete fibers, and holomorphic functions $f_{1},\cdots,f_{l}$ on X such that 
the branch locus of $\phi_{j}$
is contained in $Z_{j}=\{f_{j}=0\}$ and $Y=\displaystyle\bigcap_{j=1}^{l}Z_{j}$.\\
\\
{\bf 3. Proofs of theorems}\\
\\
\hspace*{.1in}We prove theorem 1 using the method of Fornaess and Narasimhan [8] 
(See also lemma 7, [1]).\\
\hspace*{.1in}For every irreducible component 
$X_{i}$ of X, $X_{i}\cap D\subset \subset X_{i}$ is an irreducible component
of D and a union of an increasing sequence of domains of \\
holomorphy in the Stein space $X_{i}$. (See [11]).\\
Since each $X_{i}$ is normal and $(X_{i}\cap D)_{i}$ are pairwise disjoint domains,
then we may assume that X is connected.\\
\hspace*{.1in}Let $q\in \partial{D}-Y$ and choose holomorphic functions $h_{1},\cdots,h_{m}$ on X such 
that:
$\{x\in X/h_{i}(x)=0,i=1,\cdots,m\}=\{q\}$, and j such that $q\notin 
Z_{j}$.\\
Since $D-Z_{j}$ is the union of the increasing sequence 
$(D_{k}-Z_{j})_{k\geq 1}$ of the Stein sets $D_{k}-Z_{j}$ in the Stein manifold $X-Z_{j}$,
then $D-Z_{j}$ is Stein.\newpage
\hspace*{.1in}Let $d_{j}$ be the boundary distance of the unramified domain\\
\hspace*{.2in}$\phi_{j}:D-Z_{j}\longrightarrow \complexes^{n}$. Then $-
logd_{j}$ is plurisubharmonic
on $D-Z_{j}$. Therefore the function
\hspace*{.3in}\[\psi_{j}(z)=\left\{ \begin{array}{ll}Max{(0,-
logd_{j}+k_{j}log|f_{j}|)} \ \ on \ \ D-Z_{j}\\
0 \ \ on \ \ Z_{j}
\end{array}
\right.\]
is plurisubharmonic on D, if $k_{j}$ is a large constant. This follows from a 
result due to Oka.
(See also [1], lemma 7).\\ 
\hspace*{.1in}By the Nullstellensatz, There exist a
neighborhood V of q in X and \\constants $c>0,N>0$ such that
$$\displaystyle\sum_{i=1}^{m}|h_{i}(x)|^{2}\geq c|\phi_{j}(x)-
\phi_{j}(q)|^{N},x\in V$$
\hspace*{.2in}Since $\phi_{j}$ is an analytic isomorphism at q and $f_{j}(q)\mp 0$,
it follows that,\\ if V is sufficiently small, there is a constant $c_{0}>0$ such that\\
$$\displaystyle\sum_{i=1}^{m}|h_{i}(x)|^{2}\geq c_{0}exp(-N\psi_{j}(x)), x\in 
V\cap D.$$
\hspace*{.2in}Now, since $\psi_{j}\geq 0$, there exist constants $c_{1},c_{2}>0$ 
such that\\
$$c_{2}exp(-N\psi_{j}(x))\leq \displaystyle\sum_{i=1}^{m}|h_{i}(x)|^{2}\leq 
c_{1}exp(\psi_{j}(x)), x\in D.$$
\hspace*{.2in}And applying the theorem of Skoda [14], we deduce that there is a
constant $k>0$ and holomorphic functions $g_{1},\cdots,g_{m}$ on $D-
Z_{j}$
such that
$$\displaystyle\sum_{i=1}^{m}g_{i}h_{i}=1 \ \ on \ \ D-Z_{j}$$\\

\hspace*{.7in}and\\
$$\displaystyle\sum_{i=1}^{m}\int_{D-Z_{j}}|g_{i}(x)|^{2}exp(-
k\psi_{j}(x))dv<\infty$$\newpage
Where dv is Lebesgue measure pulled back to D and k depending only
on N and m. The existence of a holomorphic function f on D  which is 
unbounded on any sequence $\{q_{\mu}\}$ of points approaching q follows from 
lemma 3-1-2 of Fornaess-Narasimhan [8]. Since $\partial{D}-Y$ is dense
in $\partial{D}$, it follows that D is a domain of holomorphy.
\\
\hspace*{.1in}We shall prove theorem 2 using the following result of 
R.Simha [15]. 
\begin{th}{-Let X be a normal Stein complex space of
dimension 2, and\\ H a hypersurface in X. Then $X-H$ is Stein.}
\end{th}
{\bf Proof of theorem 2}\\
\\
\hspace*{.1in}By the theorem of Andreotti-Narasimhan [1],
it is suffficient to prove that\\ D is locally Stein,
and we may of course assume that X is connected.\\
\hspace*{.1in}Let $p\in \partial{D}\cap Y$, and choose a connected Stein open
neighborhood U of p with $U\cap Y=\{p\}$ and such that U is biholomorphic
to a closed analytic set of a domain M in some $\complexes^{N}$.
Let E be a complex affine subspace of $\complexes^{N}$ of maximal dimension
such that p is an isolated point of $E\cap U$.\\
\hspace*{.1in}By a coordinate transformation one can obtain that $z_{i}(p)=0$
for all\\ $i\in \{1,\cdots,N\}$ and we may assume that
there is a connected Stein\\ neighborhood V of p in M such that
$U\cap V\cap \{z_{1}(x)=z_{2}(x)=0\}=\{p\}$.\\
\hspace*{.1in}We may, of course, suppose that $N\geq 4$, and let\\
$E_{1}=V\cap \{z_{1}(x)=\cdots=z_{N-2}(x)=0\}$,
$E_{2}=\{x\in E_{1}: z_{N-1}(x)=0\}$.\\
Then $A=(U\cap V)\cup E_{1}$ is a Stein closed analytic set in V
as the union of two Stein global branches of A.\\
\hspace*{.1in}Let $\zeta:\hat{A}\rightarrow A$ be a normalization of A. Then
$\zeta:\hat{A}-\zeta^{-1}(p)\rightarrow A-\{p\}$ is biholomorphic. Since
$\zeta^{-1}(E_{1})=\{x\in \hat{A}: z_{1}(\zeta(x))=\cdots=z_{N-2}(\zeta(x))=0\}$
is everywhere 1-dimensional, it follows from theorem 4 that $\hat{A}-\zeta^{-1}(E_{2})$
is Stein. Hence $A-E_{2}=\zeta(\hat{A}-\zeta^{-1}(E_{2}))$ itself is Stein.\\
\hspace*{.1in}Since $p\in E_{2}$ is the unique singular point of A, then $U\cap V\cap D$
is Stein being a domain of holomorphy in the Stein manifold $A-E_{2}$.\\
\hspace*{.1in}If $p\in \partial{D}-Y$, then there exists j such that
$p\notin Z_{j}$. We can find a Stein open neighborhood U of p in X
such that $U\cap Z_{j}=\emptyset$. Then $U\cap D=U\cap (D-Z_{j})$
is Stein.\newpage
\hspace*{.1in}The main step in the proof of theorem 3 is to show, 
when D is,\\ in addition, relatively compact in X, that for 
all $p\in \partial{D}$, there exist an open neighborhood U of p
in X and an exhaustion function f on $D\cap U$ such that for each
open $V\subset \subset U$ there is a continuous function g on V
which is locally the maximum of a finite number of stricly plurisubharmonic
functions with $|f-g|<1$. Which implies that $D\cap U$ is 1-complete with corners.\\
(A result due to Peternell [12]).\\
This result will be applied in connection with\\ the 
Diederich-Fornaess theorem [6] which asserts that an irreducible
n-dimensional complex space X is Stein if X is 1-complete with corners.\\
The proof is also based on the following 
result of M.Peternell [12]
\begin{lm}{-Let X be a complex space of pure dimension n,\\
$W\subset X\times X$ be an open set and
$f\in F_{n}(W-\Delta_{X})$ where\\ $\Delta_{X}=\{(x,x): x\in X\}$,
and let $S\subset \subset S'\subset \subset X$ be open subsets of X
such that $S\times S'\subset \subset W$.
Define $s(x)=Sup\{f(x,y): y\in \overline{S'}-S\}$ for $x\in S$ and assume
that $s(x)>f(x,y)$ if $y\in \partial{S'}$.\\
\hspace*{.1in}If S is Stein, then for each $D\subset \subset S$ and each 
$\varepsilon>0$, there is a $g\in F_{1}(D)$ such that $|g-s|<\varepsilon$
on D.\\Here $F_{n}(D)$ and $F_{1}(D)$ denote respectively the sets
of continuous functions on D which are locally the maximum of a finite
number of strongly n-convex (resp. stricly psh) functions.}
\end{lm}
{\bf Proof of theorem 3}\\
\\
\hspace*{.1in}Clearly we may suppose that D is relatively compact in X.\\
\hspace*{.1in}Since the Stein property is invariant under normalization [11],
we may\\ assume that X is normal and connected.\\
\hspace*{.1in}For $n=2$, theorem 3 follows as an immediate consequence
of theorem 2. Then we may also assume that $n\geq 3$.\\
\hspace*{.1in}Let $p\in \partial{D}\cap Y$, and choose a Stein open
neighborhood U of p in X that can be realized as a closed complex subspace
of a domain M in $\complexes^{N}$.\\
\hspace*{.1in} Let E be a complex affine subspace
of $\complexes^{N}$ of maximal dimension such that p is an isolated point of
$E\cap U$, and let E' be any complementary complex affine subspace to E in
$\complexes^{N}$ through p.\newpage
We may choose the coordinates $z_{1},\cdots,z_{N}$
and the space $E'$
such that $z_{i}(p)=0$ for all $i\in \{1,\cdots,N\}$ and $dim(E'\cap U)\geq 1$.
Since $T=E'\cap U$ is a closed analytic set in U,
and $h(z,w)=|z|^{2}+|w|^{2}-log(|z-w|^{2})$ a strongly n-convex $C^{\infty}$
function on $E'\times E'-\Delta_{E'}$, then there exists a strongly
n-convex $C^{\infty}$ function $\psi$ on a neighborhood W of 
$T'=T\times T-\Delta_{T}$ with $W\subset U\times U-\Delta_{U}$ such that
$h\leq \psi/_{T'}\leq h+1$. (See Demailly [5]).\\
\hspace*{.1in}Let W' be an open set in $U\times U$ such that $W=W'-\Delta_{U}$.
We may choose W'such that there exist a neighborhood N of p in X
and a Stein open neighborhood $U_{1}$ of p with $U_{1}\subset \subset N$
and such that $U_{1}\times (N-U_{1})\subset \subset W'$.\\
\hspace*{.1in}We now construct an exhaustion function $f_{1}$ on $U_{1}\cap D$
such that for each open $Z\subset \subset U_{1}\cap D$ there is a
$g\in F_{1}(Z)$ with $|g-f_{1}|<1$.\\
\hspace*{.1in}Let
 $f_{1}(z)=Sup\{\psi(z,w), w\in \overline{N}-U_{1}\cap D \}$,
 $z\in U_{1}\cap D$\\
\hspace*{.1in}Obviously $f_{1}$ is an exhaustion function on $U_{1}\cap D$.
There exists $m\geq 1$ such that $Z\subset \subset U_{1}\cap D_{m}$.\\
\hspace*{.1in}We now define
$$g_{j}(z)=Sup\{\psi(z,w): w\in \overline{N}-U_{1}\cap D_{j}\},
\ for \ z\in U_{1}\cap  D_{j}, \ j\geq m$$
\hspace*{.1in}Since $U_{1}\cap D_{j}$ is Stein, $\psi(z,w)$ is n-convex on W,
and $g_{j}(z)>\psi(z,w)$ for every 
$(z,w)\in (U_{1}\cap D_{j})\times \partial{N}$, then
there is a $h_{j}\in F_{1}(Z)$ such that $|g_{j}-h_{j}|<\frac{1}{2}$.
Since, obviously, $(g_{j})_{j\geq 1}$ converges uniformally on compact sets
to $f_{1}$, then there is a $j\geq m$ such that $|g_{j}-f_{1}|<\frac{1}{2}$
on Z. Hence $|f_{1}-h_{j}|<1$. Now the theorem follows from the lemma
and the theorem of Diederich-Fornaess [6]
\newpage

{\large References}\\
\\
$[1]$. A.Andreotti and R.Narasimhan. Oka's Heftungslemma and the Levi problem for 
complex spaces.Trans. AMS III(1964) 345-366.\\ 
\\
$[2]$. Behnke.H,Stein,K.: Konvergente Folgen Von Regularitatsbereichen and die 
Meromorphiekonvexitat, Math Ann. 166, 204 216(1938)\\
\\
$[3]$. M.Coltoiu,Remarques sur les r\'eunions croissantes d'ouverts de Stein 
C.R.Acad. Sci. Paris.t.307, S\'erie I,p.91-94, 1988.\\ 
\\
$[4]$. M.Coltoiu, Open problems concerning Stein spaces. Revue Roumaine de 
Math\'ematiques Pures et Appliqu\'ees.\\
\\
$[5]$. Demailly, J.P.: Cohomology of q-convex spaces in top degrees.
Math.Z 204, 283-295 (1990)\\
\\
$[6]$. Diederich, H., Fornaess, J.E.: Smoothing q-convex functions in the singular 
case. Math. Ann. 273, 665-671 (1986)\\ 
\\
$[7]$. Fornaess, J.E.: An increasing sequence of Stein manifolds whose limit is not 
Stein, Math. Ann.223, 275-277(1976)\\
\\
$[8]$. Fornaess J.E.,Narasimhan.R.:The levi problem on complex spaces with 
singularities. Math. Ann. 248(1980),47-72.\\
\\
$[9]$. Grauert, H., Remmert, R.: Singularitaten Komplexer Manngifaltigkeiten und 
Riemannsche Gebiete. Math. Z. 67, 103-128 (1957).\\ 
\\
$[10]$. Markoe, A.: Runge Families and Inductive limits of Stein spaces.Ann. Inst. 
Fourier 27,Fax.3(1977).\\
\\
$[11]$. Narasimhan, R.: A note on Stein spaces and their normalizations. Ann Scuela 
Norm. Sup. Pisa 16 (1962),327-333.
\newpage
$[12]$. Peternell, M.: Continuous q-convex exhaustion functions. Invent. Math. 85, 
246-263 (1986)\\
\\
$[13]$. Stein,K.: Uberlagerungen holomorph-vollstandiger Komplexer 
Raume.Arch.Math. 7,354-361 (1956)\\
\\
$[14]$. Skoda,H.:Application de techniques $L^{2}$ … la th\'eorie des id\'eaux d'une 
alg\'ebre de fonctions holomorphes avec poids. Ann.Sci.Ecole Norm.Sup. Paris 
5,545-579 (1972)\\
\\
$[15]$. Simha,R.: On the complement of a curve on a Stein space.
Math. Z. 82, 63-66 (1963).

\end {document}